\documentclass{article}
\usepackage{amssymb}

\usepackage{amsmath,amssymb}
\usepackage{hyperref}

\newtheorem{theorem}{Theorem}

\newtheorem{corollary}[theorem]{Corollary}

\numberwithin{equation}{section}

 \newcommand{\CP}{{\mathbb C}{\mathbb P}}
\newcommand{\R}{\mathbb R}
\begin{document}
\thispagestyle{empty}

\title{Mean curvature flows and isotopy problems}
\author{Mu-Tao Wang}

\date{April 1, 2012}
\maketitle

\begin{abstract}

In this note, we discuss the mean curvature flow of graphs of maps between Riemannian manifolds.
Special emphasis will be placed on estimates of the flow as a non-linear parabolic system of differential equations.
Several global existence theorems and applications to isotopy problems in geometry and topology will be presented.
The results are based on joint works of the author with his collaborators I. Medo\v{s}, K. Smoczyk, and M.-P. Tsui.
\end{abstract}
\section{Introduction}

We start with classical minimal surfaces in $\R^3$ (see for example \cite{os3}). Suppose a surface $\Sigma$ is given as the graph of a function $f=f(x, y)$ over a domain $\Omega\subset \R^2$ : \[\Sigma=\{(x, y, f(x, y))\,|\, (x, y)\in \Omega\}.\]  The area $A(\Sigma)$ is given by the formula \[A(\Sigma)=\int_\Omega\sqrt{1+|\nabla f|^2}.\]
The Euler-Lagrange equation for the area functional is derived to be \begin{equation}\label{mse_2} div(\frac{\nabla f}{\sqrt{1+|\nabla f|^2|}})=0.\end{equation}
Equation \eqref{mse_2}, so called the minimal surface equation, is one of the most studied nonlinear elliptic PDE and there are many beautiful classical results such as the celebrated Bernstein's conjecture for entire solutions \cite{bdm1, bdm2}. The Dirichlet problem is uniquely solvable as long as the mean curvature of the boundary $\partial \Omega$ is positive \cite{js}. In addition, any Lipschitz solution is smooth and analytic \cite{mos, mor}.

The corresponding parabolic equation is called the mean curvature flow. Here we have a time-dependent surface $\Sigma_t$,  given as the graph of a function $f=f(x, y, t)$ for each $t$, and $f$ satisfies
 { \[\frac{\partial f}{\partial t}=\sqrt{1+|\nabla f|^2} div (\frac{\nabla f}{\sqrt{1+|\nabla f|^2}}),\]} This is the negative gradient flow of the area functional. In fact, the normal component of the velocity vector of the graph of $f(x, y, t)$ in $\R^{3}$ is exactly the mean curvature vector.

The equation has been extensively studied by many authors such as Huisken \cite{hu2,hu}, Ecker-Huisken \cite{eh1, eh2}, Ilmanen \cite{il}, Andrews \cite{an}, White \cite{wh1, wh2}, Huisken-Sinestrari, \cite{hs1, hs2} X.-J. Wang \cite{xwa}, Colding-Minicozzi \cite{cm}, etc. Note that though the elliptic equation \eqref{mse_2} is in divergence form, the parabolic equation is not. Therefore, standard results from parabolic PDE theory do not readily apply.

We can also consider the equations in parametric form. Suppose the surface is given by an embedding: $\overrightarrow{X}(u, v)=(X_1(u, v), X_2(u, v), X_3(u, v))\in \R^3$.
The minimal surface equation \eqref{mse_2} is equivalent to \[(\Delta_\Sigma X_1, \Delta_\Sigma X_2, \Delta_\Sigma X_3)=(0,0,0)\] where $\Sigma$ is the image surface of $\overrightarrow{X}$  and $\Delta_\Sigma$ is the Laplace operator with respect to the induced metric on $\Sigma$. In fact, $\overrightarrow{H}=\Delta_\Sigma\overrightarrow{X}$ is the mean curvature vector of $\Sigma$. However, this elegant form has a disadvantage that it is invariant under reparametrization and thus represents a degenerate elliptic system for $(X_1, X_2, X_3)$. The same phenomenon is encountered for any curvature equation in which the diffeomorphism group appears as the symmetry group.

The corresponding parabolic equation for a family of time-dependent embeddings $\overrightarrow{X}(u, v, t)$ is
\[\frac{\partial \overrightarrow{X}}{\partial t}=\Delta_\Sigma \overrightarrow{X}.\] For this, the  mean curvature flow is often referred as the heat equation for submanifolds, just as the Ricci flow is the heat equation for Riemannian metrics. However, it is clear that the equation is of nonlinear nature as $\Delta_\Sigma$ depends on first derivatives of  $\overrightarrow{X}$.

Our subject of study in this note is a submanifold of ``higher codimension", such as a 2-surface in a 4-dimensional space given by the graph of a vector value function $(f, g)$:
\[\Sigma=\{(x, y, f(x,y), g(x, y))\,|\, (x, y)\in \Omega\}.\] The area of $\Sigma$ is then
\[A(\Sigma)=\int_\Omega \sqrt{1+|\nabla f|^2+|\nabla g|^2+{ (f_x g_y-f_yg_x)^2}}\] and the  Euler-Lagrange equation is a non-linear elliptic system for $f$ and $g$ (see the next paragraph).

In general, we consider a vector-valued function $\vec{f}:\Omega\subset \R^n\rightarrow \R^m$ and $\Sigma$ is the graph of $\vec{f}=(f^1, \cdots, f^m)$ in $\R^{n+m}$. Denote the induced metric on $\Sigma$ by \[g_{ij}=\delta_{ij}+\sum_{\alpha=1}^m \frac{\partial f^\alpha}{\partial x^i}\frac{\partial f^\alpha}{\partial x^j}.\] The volume of $\Sigma$ is \[\int_\Omega \sqrt{\det g_{ij}}\] and the Euler-Lagrange equation, which is often referred as the minimal surface system, is
\[\sum_{i, j=1}^n g^{ij} \frac{\partial^2 f^\alpha}{\partial x^i\partial x^j}=0, \,\,\alpha=1,\cdots, m,\] where $g^{ij}=(g_{ij})^{-1}$ is the inverse of $g_{ij}$.

 The corresponding parabolic equation is the mean curvature flow
\[\frac{\partial f^\alpha}{\partial t}=\sum_{i, j=1}^n g^{ij} \frac{\partial^2 f^\alpha}{\partial x^i\partial x^j}, \,\,\alpha=1,\cdots, m.\]

There is no reason to stop there and we can consider the even more general situation when ${\bf f}:M_1\rightarrow M_2$ is a differentiable map between Riemannian manifolds, and $\Sigma$ is the graph of ${\bf f}$ in $M_1\times M_2$ for $M_1$ an $n$-dimensional Riemannian manifold and $M_2$ an $m$-dimensional one.

In contrast to the codimensional one case, in an article by Lawson-Osserman \cite{lo} entitled ``Non-existence, non-uniqueness and irregularity of solutions to the minimal surface system", the undesirable features of the system mentioned in the title are demonstrated. The codimension one case, i.e. $m=1$, is essentially a scalar equation. In addition, the normal bundle of an oriented hypersurface is always trivial. On the other hand, $m>1$ corresponds to a genuine systems and the components $f^1, \cdots, f^m$ interact with each other. Moreover, the geometry of normal bundle can be rather complicated.

 Nevertheless, we managed to obtain estimates and prove global several existence theorems for higher-codimensional mean curvature flows with appropriate initial data. I shall discuss the methods in the next section before presenting the results.

 \section{Method of proofs}
Let us start with the $C^1$ estimate. In the codimension-one case (see \cite{eh1} for the equation in a slightly different but equivalent form), $m=1$, an important equation satisfied by $J_1=\frac{1}{\sqrt{1+|\nabla f|^2}}$ is  \[\frac{d}{dt}J_1=\Delta_\Sigma J_1 + R_1(\nabla f, \nabla^2 f).\] The term $R_1>0$ is quadratic in $\nabla^2 f$.

 Let us look at the $m=2, n=2$ case. We can similarly take
 \[J_2= \frac{1}{ \sqrt{1+|\nabla f|^2+|\nabla g|^2+(f_x g_y-f_yg_x)^2}}\]
 and compute the evolution equation:
\[\begin{split}\frac{d}{dt}J_2=\Delta_\Sigma J_2 + R_2(\nabla f,\nabla g, \nabla^2 f, \nabla^2 g).\end{split}\] It is observed that $R_2$ is quadratic in $\nabla^2 f$ and $\nabla^2 g$ and is positive if $|f_x g_y-f_y g_x|\leq 1$ (The is can be found in \cite{wa2}, though in a somewhat more complicated form).

A natural idea is to investigate how the quantity $f_x g_y-f_y g_x$, or the Jacobian of the map $(f, g)$ changes along the flow. Together with the maximal principle, it was shown that \cite{wa1, wa2}

(1) $f_x g_y-f_y g_x=1$ is ``preserved" along the mean curvature flow (area preserving).

(2) $|f_x g_y-f_y g_x|<1$ is ``preserved" along the mean curvature flow (area decreasing).

Here a condition is ``preserved" means if the condition holds initially, it remains true later as along as the flow exists smoothly.

 Combining with the evolution equation of $J_2$ and applying the maximum principle again show that $J_2$ has a lower bound, which in turn gives a $C^1$ estimate of $f$ and $g$. Notice that $J_2$ can be regarded as the Jacobian of the projection map onto the first factor of $\R^2$. Thus by the inverse function theorem, the graphical condition is also preserved.

 Such a condition indeed corresponds to the Gauss map of the submanifold lies in a totally geodesic or geodesically convex subset of the Grassmannian \cite{wa5}. The underlying fact for this calculation is based on the observation \cite{wa5} that the Gauss map of the mean curvature flow is a (nonlinear) harmonic map heat flow.

In codimension one case, the higher derivatives estimates follows from the $C^1$ estimates \cite{eh1}. The elliptic analogue is the theorem of Moser which states that any Lipschitz solution of the minimal surface equation is smooth. The scenario is totally different in the higher codimension case.  Lawson-Osserman constructed minimal cones in higher codimensions and thus a Lipschitz solution to the minimal surface system with $m>1$ may not be smooth at all.

Here we use ``blow-up analysis" for geometric evolution equations. An important tool is Huisken-White's monotonicity formula \cite{hu2, wh1} which characterizes central blow-up profiles as solutions of the elliptic equation:
\[\overrightarrow{H}=-\overrightarrow{X}.\]

 In general, singularity profiles for parabolic equations are soliton (self-similar) solutions of the equation. In the case of mean curvature flows, soliton (self-similar) solutions are moved by homothety or translations of the ambient space. Exclusion of self-similar ``area-preserving" or ``area-decreasing" singularity profiles and the $\epsilon$ regularity theorems of White \cite{wh2} give the desired $C^2$ estimates.

Two major difficulties remain to be overcome:

(1) Boundary value problem. This was addressed in \cite{wa6}. More sophisticated barriers that are adapted to the boundary geometry are needed in order to obtain sharper result to cover the area-decreasing case.

(2) Effective estimates in time as $t\rightarrow \infty$. So far, convergence results rely on the sign of the curvature of the ambient space. The $C^2$ estimates obtained through blow-up analysis usually deteriorate in time.

In the next section, we present the statements of results which are cleanest when $M_1$ and $M_2$ are closed Riemannian manifolds with suitable curvature conditions. We remark that there have been several global existence and convergence theorems on higher codimensional graphical mean curvature flows such as \cite{sm2,sw, wa3, clt, ab}, etc. Here we focus on those theorems that have implications on isotopy problems.

\section{Statements of results related to isotopy problems}

\subsection{Symplectomorphisms of Riemann surfaces}

Let $(M_1, g_1)$ and $(M_2, g_2)$ be Riemann surfaces with metrics of the same constant curvature. We can normalize so the curvature is $-1, 0$ or $1$.
Let $f:M_1\rightarrow M_2$ be an oriented area-preserving map and $\Sigma$ be the graph of $f$ in $M_1\times M_2 $. A oriented area-preserving map is also a symplectomorphism, i.e. $f^*\omega_2=\omega_1$ where $\omega_1$ and $\omega_2$ are the area forms (or symplectic forms) of $g_1$ and $g_2$, respectively. The area $A(f)$ of the graph of $f$ is a symmetric function on the symplectomorphism group, i.e. $A(f)=A(f^{-1})$ and the mean curvature flow gives a deformation retract of this group to a finite dimensional one.

\begin{theorem}{(\cite{wa1, wa2, wa4}, see also \cite{wa7}) Suppose $\Sigma_0$ is the graph of a symplectomorphism $f_0:M_1\rightarrow M_2$. The mean curvature flow $\Sigma_t$ exists for all $t\in [0, \infty)$ and converges smoothly to a minimal submanifold as $t\rightarrow \infty$.
$\Sigma_t$ is the graph of  a symplectic isotopy $f_t$ from $f_0$ to a canonical minimal map $f_\infty$.}
\end{theorem}

Since any diffeomorphism is isotopic to an area preserving diffeomorphism,
this gives a new proof of Smale's theorem \cite{smale} that O(3) is the deformation
retract of the diffeomorphism group of $S^2$. For a positive genus Riemann surface,
this implies the identity component of the diffeomorphism group is contractible.

The result for the positive genus case was also obtained by Smoczyk \cite{sm2} under an extra angle condition.

In this case, the graph of the symplectomorphism is indeed a Lagrangian submanifolds in the product space. There have been important recent progresses on the Lagrangian minimal surface equation, we refer to the excellent survey article of Brendle \cite{br} in this direction.

For an area-decreasing map $f$, i.e. $|f^*\omega_2|<\omega_1$, the mean curvature flow exists for all time and converges to the graph of a constant map, see \cite{wa2}.

\subsection{Area-decreasing maps in higher dimensions}
 The area-decreasing condition, which turns out to be rather natural for the mean curvature flow, can be generalized to higher dimensions. A Lipschitz map $f:M_1\rightarrow M_2$ between Riemannian manifolds is area-decreasing if the 2-dilation $|\Lambda^2 df|_p|<1 $ for each $p\in M_1$. Here $\Lambda^2 df|_p:\Lambda^2 T_pM_1\rightarrow \Lambda^2 T_{f(p)} M_2$ is the map induced by the differential $df|_p: T_pM_1\rightarrow T_{f(p)} M_2$.

 Equivalently, in local orthonormal coordinate systems on the domain and the target, we ask
    \[|\frac{\partial f^\alpha}{\partial x^i}\frac{\partial f^\beta}{\partial x^j}-\frac{\partial f^\alpha}{\partial x^i}\frac{\partial f^\beta}{\partial x^j}|<1\] for $\alpha\not=\beta$, $i\not= j$.  This is also the same as  ${H}^2(f(D))\leq {H}^2(D)$ for any $D\subset M_1$ of finite two-dimensional Hausdorff measure $H^2(\cdot)$.

In \cite{tw}, we proved that area decreasing condition is preserved along the mean curvature flow for the graph of a smooth map $f:S^n\rightarrow S^m$ between spheres of constant curvature $1$. In addition,

\begin{theorem}\cite{tw} Suppose $n, m\geq 2$. If $f:S^n\rightarrow S^m$ is an area-deceasing smooth map, the mean curvature flow of the graph of $f$ exists for all time, remains a graph, and converges smoothly to a constant map as $t\rightarrow \infty$.\end{theorem}

 The most difficult part of the proof is to express the area-decreasing condition as the two-positivity condition (i.e. the sum of the two smallest eigenvalues is positive) for a Lorentzian metric of signature $(n, m)$ and compute the evolution equation of the induced metric.

A simple corollary is the following:

\begin{corollary}  If $n, m\geq 2$, every area-decreasing map $f:S^n\rightarrow S^m$ is homotopically trivial. \end{corollary}

Gromov \cite{gr2} shows that for each pair $(n, m)$, there exists a number
$\epsilon(n,m) > 0$, so that any map from $S^n$ to $S^m$
with $|\Lambda^2 df| < \epsilon(n,m)$ is homotopically trivial, where $\epsilon (n, m)<< 1$. In general, we may consider the $k$-Jacobian $\Lambda^k d f:\Lambda^k TM_1\rightarrow \Lambda^k TM_2$, whose supreme norm $|\Lambda^k df|$ is called the $k$-dilation ($k=1$ is the Lipschitz norm). Guth \cite{lg} constructed homotopically non-trivial maps from $S^n$ to $S^m$ with arbitrarily small 3-dilation. It is amazing that 2-dilation is sharp here as it arises naturally from a completely different consideration of the Gauss map of the mean curvature flow (see last section).

\subsection{Symplectomorphisms of complex projective spaces}
In this section, we consider the generalization of the theorem for symplectomorphisms of Riemann surfaces to higher dimensional manifolds.
Let $M_1$ and $M_2$ be K\"ahler manifolds equipped with K\"ahler-Einstein metrics of the same Ricci curvature. Let $f:M_1\rightarrow M_2$ be a symplectomorphism. As was remarked in the last section, we can consider the graph of $f$ as a Lagrangian submanifold $\Sigma$ in the product space $M_1\times M_2$ and deform it by the mean curvature flow.  A theorem of Smoczyk \cite{sm1} (see also \cite{oh}) implies that the mean curvature flow $\Sigma_t$ remains a Lagrangian submanifold. If we can show $\Sigma_t$ remains graphical as well, it will corresponds to a symplectic isotopy $f_t:M_1\rightarrow M_2$.  The simplest case to be considered in higher dimension is  $M_1=M_2=\CP^n$ with the Fubini-Study metric. In a joint work with Medo\v{s}, we proved the following pinching theorem.

\begin{theorem}\cite{mw} There exists an explicitly computable constant $\Lambda>1$ depending only on $n$, such that any symplectomorphism $f:\CP^n\rightarrow \CP^n$ with \[\frac{1}{\Lambda} g\leq f^* g\leq \Lambda g\]
is symplectically isotopic to a biholomorphic isometry of $\CP^n$ through the mean curvature flow.\end{theorem}

A theorem of Gromov \cite{gr1} shows that, when $n=2$, the statement holds true without any pinching condition by the method of pseudoholomorphic curves. Our theorem is not strong enough to give an analytic proof of Gromov's theorem for $n=2$. However, for $n\geq 3$, this seems to be the first known result.

  Unlike previous theorems, Grassmannian geometry does not quite help here, as the subset  that corresponds to biholomorphic isometries does not have any convex neighborhood in the Grassmannian. The integrability condition, or  the Gauss-Codazzi equations, is used in an essential way to overcome this difficulty.


\end{document}